\documentclass{article}%
\usepackage{amssymb}
\usepackage{amsmath}
\usepackage{amsfonts}
\usepackage{graphicx}%
\providecommand{\U}[1]{\protect\rule{.1in}{.1in}}
\newtheorem{theorem}{Theorem}

\newtheorem{corollary}[theorem]{Corollary}

\newtheorem{lemma}[theorem]{Lemma}

\newtheorem{problem}[theorem]{Problem}

\newtheorem{remark}[theorem]{Remark}

\begin{document}
To appear in Electronic Communications in Probability\bigskip\ (2020)

\begin{center}
{\large On the CLT for additive functionals of Markov chains }

\bigskip

Magda Peligrad

\bigskip

Department of Mathematical Sciences, University of Cincinnati, PO Box 210025,
Cincinnati, Oh 45221-0025, USA. \texttt{ }
\end{center}

email: peligrm@ucmail.uc.edu

\bigskip


\smallskip

\bigskip

\begin{center}

Abstract
\end{center}

In this paper we study the additive functionals of Markov chains via
conditioning with respect to both past and future of the chain. We shall point
out new sufficient projective conditions, which assure that the variance of
partial sums of $n$ consecutive random variables of a stationary Markov chain
is linear in $n$. The paper also addresses the central limit theorem problem
and is listing several open questions.

\section{Introduction}

Throughout the paper assume that $(\xi_{n})_{n\in\mathbb{Z}}$ is a stationary
Markov chain defined on a probability space $(\Omega,\mathcal{F},\mathbb{P})$
with values in a measurable space $(S,\mathcal{A})$. We suppose that there is
a regular conditional distribution for $\xi_{1}$ given $\xi_{0}$ denoted by
$Q(x,A)=\mathbb{P}(\xi_{1}\in A|\,\xi_{0}=x)$. In addition $Q$ denotes the
Markov operator {acting via $(Qf)(x)=\int_{S}f(s)Q(x,ds)$. }Denote by
$\mathcal{F}_{n}=\sigma(\xi_{k},k\leq n)$ and $\mathcal{F}^{n}=\sigma(\xi
_{k},k\geq n)$. The invariant distribution is denoted by $\pi(A)=\mathbb{P}%
(\xi_{0}\in A)$ and { }$Q^{\ast}$ denotes the adjoint of $Q$.{ Next, let
$\mathbb{L}_{0}^{2}(\pi)$ be the set of measurable functions on $S$ such that
$\int f^{2}d\pi<\infty$ and $\int fd\pi=0.$ For a function} ${f}\in
${$\mathbb{L}_{0}^{2}(\pi)$ let }%
\[
{X_{i}=f(\xi_{i}),\ S_{n}=\sum\limits_{i=1}^{n}X_{i}.}%
\]
{ Note that for every $k\in\mathbb{Z},$ $Q^{k}f(\xi_{0})=\mathbb{E}(X_{k}%
|\xi_{0}),$ while $(Q^{\ast})^{k}f(\xi_{0})=\mathbb{E}(X_{-k}|\xi_{0})$. We
denote by }${{||X||}}$ the norm in {$\mathbb{L}^{2}$}$(\mathbb{P})$ and by
$||f||_{\pi}$ the norm in {$\mathbb{L}_{0}^{2}(\pi)$}.{ Sometimes, we shall
also use the notations }%
\[
V_{n}=I+Q+...+Q^{n}\,\text{\ and }V_{n}^{\ast}=I+Q^{\ast}+...+(Q^{\ast})^{n}.
\]
In some statements, we assume that the stationary Markov chain is ergodic,
i.e. the only invariant functions $Qf=f$ are the constant functions. {
Concerning the central limit theorem for additive functionals of a stationary
and ergodic Markov chain, many of the results in the literature are given
under sufficient conditions either in terms of }$Q^{k}f$ or $V_{n}f.$ Among
them we mention the pioneering works by Gordin (1969), Gordin and Lifshitz
(1978), Heyde (1974), McLeish (1975) and Voln\'{y} (1993) among others. For a
survey see Peligrad (2010) and the book by Merlev\`{e}de et al. (2019).

Maxwell and Woodroofe (2000) introduced a more general condition than in the
papers mentioned above, namely%
\begin{equation}
\sum_{n\geq1}\frac{||E(S_{n}|\xi_{0})||}{n^{3/2}}<\infty. \label{MW}%
\end{equation}
In the same paper, they showed that (\ref{MW}) implies that%
\begin{equation}
\frac{E(S_{n}^{2})}{n}\rightarrow\sigma^{2} \label{varSn}%
\end{equation}
and
\begin{equation}
\frac{S_{n}}{\sqrt{n}}\Rightarrow N(0,\sigma^{2}), \label{CLT}%
\end{equation}
where $\Rightarrow$ denotes the convergence in distribution and $N(0,\sigma
^{2})$ is a normally distributed random variable. Later on, Peligrad and Utev
(2005)\ established the functional form of the CLT under (\ref{MW}).

There are examples of Markov chains pointing out that, in general, condition
(\ref{MW}) is as sharp as possible is some sense. Peligrad and Utev (2005)
constructed an example showing that for any sequence of positive constants
$(a_{n}),$ $a_{n}\rightarrow0,$ there exists a stationary Markov chain such
that
\[
\sum_{n\geq1}a_{n}\frac{||E(S_{n}|\xi_{0})||}{n^{3/2}}<\infty
\]
but $S_{n}/\sqrt{n}$ is not stochastically bounded. This example and other
counterexamples provided by Voln\'{y} (2010), Dedecker (2015) and Cuny and Lin
(2016), show that, in general, condition
\begin{equation}
\sum_{n\geq1}\frac{||E(S_{n}|\xi_{0})||^{2}}{n^{2}}<\infty\label{1}%
\end{equation}
does not assure that (\ref{varSn}) holds and also does not assure (\ref{CLT}).

However, by using Proposition 2.2 in Cuny (2011), which connects
(\ref{1})\ with a spectral condition, we know that (\ref{1}) is sufficient for
the CLT given in (\ref{CLT}),\ in case when the Markov chain is normal
($QQ^{\ast}=Q^{\ast}Q$), as announced by Gordin and Lifshitz (1981) and proven
in Section IV.7 in Borodin and Ibragimov (1995)\ and also, independently, in
Derriennic and Lin (1996).

\begin{theorem}
\label{CLTnormal}Gordin and Lifshitz (1981). Assume that the Markov chain is
normal, stationary and ergodic and satisfies (\ref{1}). Then (\ref{varSn}) and
(\ref{CLT}) hold.
\end{theorem}

For reversible Markov chains ($Q=Q^{\ast}$) condition (\ref{1})\ is equivalent
to (\ref{varSn}) and also with the convergence of $E_{\pi}(fV_{n}f)$ (see
Kipnis and Varadhan (1986), Proposition 2.2 in Cuny (2011) and the remarks
following Proposition 1 in Derriennic and Lin \cite{DL2}). Furthermore in this
case the functional CLT\ also holds.

A natural question is to ask what will be a natural generalization of Theorem
\ref{CLTnormal} to Markov processes, which are not necessarily normal. In
other words what will be a natural minimal condition to be added to (\ref{1}),
which will insure (\ref{varSn}) and (\ref{CLT}).

A possibility is to impose besides (\ref{1}) a similar condition, but
conditioning this time with respect to the future of the process
\begin{equation}
\sum_{n\geq1}\frac{||E(S_{n}|\xi_{n})||^{2}}{n^{2}}<\infty. \label{2}%
\end{equation}
{ }It is interesting to point out that for normal Markov chains conditions
(\ref{1}) and (\ref{2}) coincide.

In the operator notation, conditions (\ref{1}) and (\ref{2}) could be written
in the following alternative form:%
\begin{equation}
\sum_{n\geq1}\frac{||V_{n}(f)||_{\pi}^{2}}{n^{2}}<\infty\text{ and }%
\sum_{n\geq1}\frac{||V_{n}^{\ast}(f)||_{\pi}^{2}}{n^{2}}<\infty.
\label{two cond}%
\end{equation}
This paper has double scope. First, in Section 2, we shall raise some open
questions concerning the CLT for the additive functionals of a Markov chain
under conditions related to (\ref{two cond}). In the following section we
support these conjectures by proving some partial results. We shall show, for
instance, that (\ref{two cond}) implies (\ref{varSn})\ and we shall comment
that the CLT holds up to a random centering. The proofs are given in Sections
4 and 5.

\section{Open problems}

We shall list here several natural open problems.

\begin{problem}
\label{Pr1}For a stationary and ergodic Markov chain is it true (or not) that
condition (\ref{two cond}) implies that the CLT in (\ref{CLT}) holds?
\end{problem}

In terms of the individual random variables, let us note that, by the triangle
inequality, stationarity and Lemma \ref{aux} in the last section, applied with
$a_{k}=||E(X_{k}|\xi_{0})||$ (and also with $a_{k}=||E(X_{-k}|\xi_{0})||^{2}%
$), we obtain%

\[
\sum\nolimits_{n\geq1}\frac{||E(S_{n}|\xi_{0})||^{2}}{n^{2}}\leq
4\sum\nolimits_{k\geq1}||E(X_{k}|\xi_{0})||^{2}%
\]
and also%
\[
\sum\nolimits_{n\geq1}\frac{||E(S_{n}|\xi_{n})||^{2}}{n^{2}}\leq
4\sum\nolimits_{k\geq1}||E(X_{-k}|\xi_{0})||^{2}.
\]
Then, clearly (\ref{1}) is implied by%
\begin{equation}
\sum\nolimits_{k\geq1}||E(X_{k}|\xi_{0})||^{2}<\infty\label{mix1}%
\end{equation}
and (\ref{2}) is implied by
\begin{equation}
\sum\nolimits_{k\geq1}||E(X_{-k}|\xi_{0})||^{2}<\infty. \label{mix2}%
\end{equation}
These two last conditions can be reformulated as:
\begin{equation}
\sum\nolimits_{k\geq1}||Q^{k}f||_{\pi}^{2}<\infty\text{ and }\sum
\nolimits_{k\geq1}||(Q^{\ast})^{k}f||_{\pi}^{2}<\infty. \label{two mix}%
\end{equation}
As a matter of fact there are summability conditions that interpolates between
(\ref{two cond})\ and (\ref{two mix}), which are known under the name of
square root conditions. Following Derriennic and Lin \cite{DL} the operator
$\sqrt{I-Q}$ is defined by%
\[
\sqrt{I-Q}:=I-\sum\nolimits_{n\geq1}\delta_{n}Q^{n},
\]
where $\sqrt{1-x}=1-\sum\nolimits_{n\geq1}\delta_{n}x^{n},$ with $\delta
_{n}>0,$ $n\geq1$ and $\sum\nolimits_{n\geq1}\delta_{n}=1.$ By the equivalent
definitions in Corollary 2.12 in Derriennic and Lin \cite{DL}, $f\in\sqrt
{1-Q}\mathbb{L}_{2}(\pi)$ is equivalent to
\[
\sum\nolimits_{k=1}^{n}\frac{1}{k^{1/2}}Q^{k}f\ \text{converges in }%
\mathbb{L}_{2}(\pi).
\]

By Proposition 4.6 of Cohen et al. (2017), applied with $b_{n}=(1+n)^{-1/2}$,
we know that $f\in\sqrt{1-Q}\mathbb{L}_{2}(\pi)$ implies (\ref{1}) and also
$f\in\sqrt{1-Q^{\ast}}\mathbb{L}_{2}(\pi)$ implies (\ref{2}). On the other
hand, by Proposition 9.2 in Cuny and Lin (2016), condition (\ref{mix1})
implies $f\in\sqrt{1-Q}\mathbb{L}_{2}(\pi)$ and condition (\ref{mix2}) implies
$f\in\sqrt{1-Q^{\ast}}\mathbb{L}_{2}(\pi).$ Furthermore, according to
Corollary 4.7 in Cohen et al. (2017), for normal contractions $f\in\sqrt
{1-Q}\mathbb{L}_{2}(\pi)$ is equivalent to (\ref{1}). We mention that
Voln\'{y} (2010) constructed an example of a (non-normal) Markov operator $Q$
and $f\in\sqrt{1-Q}\mathbb{L}_{2}(\pi)$ for which the asymptotic variance of
$||S_{n}||^{2}/n$\ does not exist. Note however that if $f\in\sqrt
{1-Q}\mathbb{L}_{2}(\pi)\cap\sqrt{1-Q^{\ast}}\mathbb{L}_{2}(\pi)$ then
(\ref{varSn}) holds (see Proposition 1 in \cite{DL2} and the remarks following
this proposition).

These considerations suggest that the following conjecture deserves to be
studied, of course, in case the answer to Problem \ref{Pr1} is negative.

\begin{problem}
\label{Pr2}If the Markov chain is stationary and ergodic is it true (or not)
that $f\in\sqrt{1-Q}\mathbb{L}_{2}(\pi)\cap\sqrt{1-Q^{\ast}}\mathbb{L}_{2}%
(\pi)$ implies that the CLT\ in (\ref{CLT}) holds?
\end{problem}

Finally, if the answer to Problem \ref{Pr2} is negative we could ask the
following question:

\begin{problem}
\label{Pr3}If the stationary Markov chain is ergodic in the ergodic
theoretical sense is it true (or not) that condition (\ref{two mix}) implies
that (\ref{CLT}) holds?
\end{problem}

\section{Results}

We give here a few results in support of the open problems which have been
raised in the previous section. Point (a) of the next theorem deals with the
variance of partial sums, which plays a very important role in the CLT. In the
next theorem by totally ergodic chain we understand that $Q^{k}$ is ergodic
for any$\ $positive integer $k$.

\begin{theorem}
\label{Thvar}Assume that conditions (\ref{1}) and (\ref{2}) hold.
Then:\newline(a) The limit in (\ref{varSn}) holds, namely\newline%
\[
\lim_{n\rightarrow\infty}\frac{E(S_{n}^{2})}{n}=\sigma^{2}.
\]
(b) If the chain is totally ergodic then the following limit exists%
\begin{equation}
\lim_{n\rightarrow\infty}\frac{1}{n}||E(S_{n}|\xi_{0},\xi_{n})||^{2}%
=\mathbb{\eta}^{2} \label{sigma-eta1}%
\end{equation}
and\newline(c)%
\begin{equation}
\frac{S_{n}-E(S_{n}|\xi_{0},\xi_{n})}{\sqrt{n}}\Rightarrow N(0,\sigma
^{2}-\mathbb{\eta}^{2}). \label{T Centered CLT}%
\end{equation}

\end{theorem}

As an immediate consequence, by the discussion in the previous section, we
also have the following corollaries:

\begin{corollary}
\label{corvar1}The conclusion of Theorem \ref{Thvar} also holds for $f\in
\sqrt{1-Q}\mathbb{L}_{2}(\pi)\cap\sqrt{1-Q^{\ast}}\mathbb{L}_{2}(\pi).$
\end{corollary}

\begin{corollary}
\label{corvar2}The conclusion of Theorem \ref{Thvar} also holds under the
couple of conditions (\ref{mix1}) and (\ref{mix2}).
\end{corollary}

\begin{remark}
\label{rem1}Let us mention that the conclusion of Corollary \ref{corvar1} does
not hold if we assume only that $f\in\sqrt{1-Q}\mathbb{L}_{2}(\pi)$ as shown
in Voln\'{y} (2010). Also, the conclusion of Corollary \ref{corvar2} does not
hold under just (\ref{mix1}). Dedecker (2015) constructed a relevant example,
which has been reformulated in Proposition 9.5. in Cuny and Lin (2016). This
example shows that there exists a Markov operator $Q$ on some $\mathbb{L}%
_{2}(\pi)$ and a function $f\in\mathbb{L}_{0}^{2}(\pi)$ satisfying
\[
\sum\nolimits_{k\geq1}(\log k)||Q^{k}f||_{\pi}^{2}<\infty
\]
and such that $||S_{n}||^{2}/n\rightarrow\infty$ as $n\rightarrow\infty.$
\end{remark}

An interesting question asked in Problem \ref{Pr1} is whether the random
centering in the point (c) of Theorem (\ref{Thvar}) can be avoided altogether.
We can prove this fact under the condition
\begin{equation}
\sum_{n\geq1}\frac{||E(S_{n}|\xi_{0},\xi_{n})||^{2}}{n^{2}}<\infty,
\label{bad}%
\end{equation}
namely:

\begin{corollary}
\label{CorBad}Assume (\ref{bad}) holds. Then (\ref{varSn}) and (\ref{CLT}) hold.
\end{corollary}

By Lemma \ref{aux}, relation (\ref{bad}) is implied by%

\begin{equation}
\sum\nolimits_{k\geq1}||E(X_{0}|\xi_{-k},\xi_{k}))||^{2}<\infty.
\label{mixingale}%
\end{equation}
Therefore we obtain the following corollary which was also pointed out in
Peligrad (2020):

\begin{corollary}
\label{cor mixi3}Assume that (\ref{mixingale}) holds.\ Then (\ref{varSn}) and
(\ref{CLT}) hold.
\end{corollary}

\section{Proofs}

Let us comment first about conditions (\ref{1}) and (\ref{2}). We are going to
establish two lemmas (Lemma \ref{Lnegli} and Lemma \ref{Lsubad}) showing that
condition (\ref{1}) implies that
\begin{equation}
\sum_{i\geq0}\frac{||E(S_{2^{i}}|\xi_{0})||^{2}}{2^{i}}<\infty, \label{2prime}%
\end{equation}
while condition (\ref{2}) implies that%
\begin{equation}
\sum_{i\geq0}\frac{||E(S_{2^{i}}|\xi_{2^{i}})||^{2}}{2^{i}}<\infty.
\label{3prime}%
\end{equation}
As a matter of fact, we can replace conditions (\ref{1}) and (\ref{2}) in
Theorem \ref{Thvar} by conditions (\ref{2prime}) and (\ref{3prime}).

\bigskip

\textbf{Proof of point (a) of Theorem \ref{Thvar}\ }

\bigskip

The proof of point (a) of Theorem \ref{Thvar}\ is related to the proof of
Proposition 2.1 in Peligrad and Utev (2005), but it takes advantage of the
Markov property. It includes several steps.

\bigskip

\textbf{1.} \textbf{Upper bound on a subsequence}

\bigskip

We shall establish first the following recurrence formula which has interest
in itself.

Denote
\[
\Delta_{2^{r}}=\sum_{i=0}^{r-1}\frac{||E(S_{2^{i}}|\xi_{0})||\cdot
||E(S_{2^{i}}|\xi_{2^{i}})||}{2^{i}}.
\]
Then, for $2^{r-1}\leq n<2^{r},$ we have the following\ bound:%
\begin{equation}
E(S_{2^{r}}^{2})\leq2^{r}[E(X_{0}^{2})+\Delta_{2^{r}}]. \label{diatic}%
\end{equation}
To establish it, denote by $\bar{S}_{n}=\sum_{k=n+1}^{2n}\xi_{k}.$ So, by
stationarity%
\[
E(S_{2n}^{2})=2E(S_{n}^{2})+2E(S_{n}\bar{S}_{n}).
\]
Note that, by the properties of conditional expectation and by the Markov
property,
\[
E(S_{n}\bar{S}_{n})=E[E(S_{n}E(\bar{S}_{n}|\mathcal{F}_{n})]=E[E(S_{n}%
E(\bar{S}_{n}|\mathcal{\xi}_{n})]=E[E(S_{n}|\xi_{n})E(\bar{S}_{n}|\xi_{n})].
\]
We see that, by recurrence we have
\begin{equation}
E(S_{2^{r}}^{2})=2^{r}\left[  E(X_{0}^{2})+\sum\nolimits_{j=0}^{r-1}\frac
{1}{2^{j}}E[E(S_{2^{j}}|\xi_{2^{j}})E(\bar{S}_{2^{j}}|\xi_{2^{j}})]\right]  .
\label{rec}%
\end{equation}
Now, by H\"{o}lder's inequality and stationarity,%
\[
|E[E(S_{n}|\xi_{n})E(\bar{S}_{n}|\xi_{n})]|\leq||E(S_{n}|\xi_{n}%
)||\cdot||E(S_{n}|\xi_{0})||,
\]
and (\ref{diatic}) follows.

\bigskip

2. \textbf{Limit on a subsequence}

\bigskip

Note that, if $\sup_{r}\Delta_{2^{r}}<\infty,$ then $\sum\nolimits_{j=0}%
^{r-1}2^{-j}E[E(S_{2^{j}}|\xi_{0})E(\bar{S}_{2^{j}}|\xi_{2^{j}})]$ converges
as $r\rightarrow\infty$, say to $L$. Then, by (\ref{rec}), we have that
\begin{equation}
\frac{1}{2^{r}}E(S_{2^{r}}^{2})\rightarrow\sigma^{2}\text{ as }r\rightarrow
\infty, \label{conv1}%
\end{equation}
where $\sigma^{2}=E(X_{0}^{2})+L$\textbf{.}

\bigskip

3.\textbf{ Limiting variance for }$S_{n}/\sqrt{n}$

\bigskip

We show here that if conditions (\ref{1}) and (\ref{2}) hold then
(\ref{varSn}) holds.

For $2^{r-1}\leq n<2^{r}$, we use the binary expansion
\begin{equation}
n=\sum_{k=0}^{r-1}2^{k}a_{k}\quad\text{where}\quad a_{r-1}=1\quad
\text{and}\quad a_{k}\in\{0,1\}. \label{repn}%
\end{equation}
Then, we apply the following representation
\begin{equation}
S_{n}=\sum_{j=0}^{r-1}U_{2^{j}}a_{j}\quad\text{where}\quad U_{2^{j}}%
=\sum_{i=n_{j-1}+1}^{n_{j}}X_{i}\;,\;n_{j}=\sum_{k=0}^{j}2^{k}a_{k}\;,\text{
}n_{-1}=0\;. \label{repr}%
\end{equation}
Clearly, for $a_{j}=0$, $U_{2^{j}}=0.$

Then we use the representation
\begin{equation}
E(S_{n}^{2})=\sum_{i=0}^{r-1}a_{i}E(U_{2^{i}}^{2})+\sum_{i\neq j=0}^{r-1}%
a_{i}a_{j}E(U_{2^{i}}U_{2^{j}})\equiv I_{n}+J_{n}. \label{estS_n}%
\end{equation}
Now, by stationarity, the representation of $n$ in (\ref{repn})\ and the
convergence in (\ref{conv1}) we obtain
\[
{\frac{I_{n}}{n}=}\frac{1}{n}\sum_{i=0}^{r-1}a_{i}2^{i}\left(  \frac
{E(S_{2^{i}}^{2})}{2^{i}}-\sigma^{2}\right)  +\sigma^{2}\rightarrow\sigma
^{2}\text{ as }n\rightarrow\infty\mathbf{.}%
\]
It remains to prove that$\mathbf{\quad|}J_{n}|/n\rightarrow0.$ Let $0\leq
i<j<r$. Then, by the properties of Markov chains and H\"{o}lder's inequality
\begin{equation}
|E(U_{2^{i}}U_{2^{j}})|\leq|E(U_{2^{i}}|\xi_{n_{i}})E(U_{2^{j}}|\xi_{n_{j-1}%
})|\leq\Vert E(S_{2^{i}}|\xi_{2^{i}})\Vert\cdot\Vert E(S_{2^{j}}|\xi_{0}%
)\Vert. \label{Holder}%
\end{equation}
Hence,%
\begin{align*}
|J_{n}|  &  \leq\frac{1}{2}\sum_{0\leq i\neq j\leq r-1}2^{i/2}2^{j/2}\left(
\frac{\Vert E(S_{2^{i}}|\xi_{2^{i}})\Vert^{2}}{2^{i}}+\frac{\Vert E(S_{2^{j}%
}|\xi_{0})\Vert^{2}}{2^{j}}\right)  \leq\\
&  \frac{1}{2(\sqrt{2}-1)}2^{r/2}\left(  \sum_{i=0}^{r-1}2^{i/2}\frac{\Vert
E(S_{2^{i}}|\xi_{2^{i}})\Vert^{2}}{2^{i}}+\sum_{j=0}^{r-1}2^{j/2}\frac{\Vert
E(S_{2^{j}}|\xi_{0})\Vert^{2}}{2^{j}}\right)  .
\end{align*}
We can easily see that $E|J_{n}|/n\rightarrow0$ because of (\ref{2prime}) and
(\ref{3prime}). $\ \square$

\bigskip\ 

\textbf{Proof of points (b) and (c) of Theorem \ref{Thvar}\ }

\bigskip

For obtaining the points (b) and (c) of Theorem \ref{Thvar}, we shall combine
the result in point (a) with a new CLT for additive functionals of Markov
chains given by Peligrad (2020), namely:

\begin{theorem}
\label{Th random center}Peligrad (2020). Assume that if the chain is totally
ergodic and
\begin{equation}
\sup_{n\geq1}\frac{E(S_{n}^{2})}{n}<\infty. \label{varsup1}%
\end{equation}
Then, the following limit exists
\begin{equation}
\lim_{n\rightarrow\infty}\frac{1}{n}||S_{n}-E(S_{n}|\xi_{0},\xi_{n}%
)||^{2}=\mathbb{\theta}^{2} \label{def teta}%
\end{equation}
and
\[
\frac{S_{n}-E(S_{n}|\xi_{0},\xi_{n})}{\sqrt{n}}\Rightarrow N(0,\mathbb{\theta
}^{2}).
\]

\end{theorem}

Note that, since (\ref{varSn}) holds, then clearly (\ref{varsup1}) holds and
in addition (\ref{def teta}) implies that
\[
\lim_{n\rightarrow\infty}\frac{1}{n}||E(S_{n}|\xi_{0},\xi_{n})||^{2}%
=\sigma^{2}-\mathbb{\theta}^{2}.
\]
It follows that points (b) and (c) of Theorem \ref{Thvar} follow with
$\eta^{2}=\sigma^{2}-\mathbb{\theta}^{2}.$

\bigskip

\textbf{Proof of Corollary \ref{CorBad}}.

\bigskip

\ By the properties of the conditional expectation we see that condition
(\ref{bad}) implies that (\ref{1}) and (\ref{2}) are satisfied and therefore,
by point (b) of Theorem \ref{Thvar}, the limit in (\ref{sigma-eta1})\ exists.
If this limit is not $0,$ note that (\ref{bad}) cannot be satisfied. Therefore
(\ref{bad}) implies that the limit in (\ref{sigma-eta1})\ is $0.$ We can apply
now Theorem 3.1 in Billingsley (1999) to conclude that
in\ (\ref{T Centered CLT})\ the random centering is not needed if we assume
(\ref{bad}).

\section{Auxiliary results}

The following lemma holds for any subadditive sequence $(V_{m})_{m\geq1}$ of
positive numbers. Its proof is inspired by Lemma 2.8. in Peligrad and Utev
(2005). Because of the subtle differences\ we shall give it here. The main
difference is that the sequence $V_{m}^{2}$ is not subadditive.

\begin{lemma}
\bigskip\ \label{Lnegli}For any positive subadditive sequence $(V_{m}%
)_{m\geq1}$ of positive numbers we have
\[
\sum_{i\geq1}\frac{V_{2^{i}}^{2}}{2^{i}}\leq65\sum_{k\geq1}\frac{V_{k}^{2}%
}{k^{2}}.
\]

\end{lemma}

Proof. We recall first a property on the page 806 in Peligrad and Utev (2005).
For a finite set of real numbers $C,$ denote by $|C|$ its cardinal. Consider a
positive integer $N$ and the set
\[
A_{N}=\{1\leq i\leq N\;:\;V_{i}\geq V_{N}/2\}.\;
\]

\textit{Property \/}: $|A_{N}|\geq N/2;$ that is $A_{N}$ contains at least
$N/2$ elements.

We start the proof of this lemma by adding the variables in blocks in the
following way:
\begin{equation}
\sum_{n\geq2}\frac{V_{n}^{2}}{n^{2}}=\sum_{r\geq0}\sum_{n=4^{r}+1}^{4^{r+1}%
}\frac{V_{n}^{2}}{n^{2}}\geq\sum_{r\geq0}\frac{1}{4^{2r+2}}\sum_{n=4^{r}%
+1}^{4^{r+1}}V_{n}^{2}. \label{5}%
\end{equation}
Define
\[
C_{r}=\{n\in\{4^{r}+1,\ldots,4^{r+1}\}\;:\;V_{n}\geq V_{4^{r+1}}%
/2\}=A_{4^{r+1}}\cap\{4^{r}+1,\ldots,4^{r+1}\}.
\]
Note that, by applying the above property with $N=4^{r+1},$ it is easy to see
that
\[
|C_{r}|\geq|A_{4^{r+1}}|-|\{1,2,...,4^{r}\}|\geq4^{r+1}/2-4^{r}=4^{r}.
\]
It follows that
\[
\sum_{n\geq2}\frac{V_{n}^{2}}{n^{2}}\geq\sum_{r\geq0}\frac{4^{r}}{4^{2r+4}%
}V_{4^{r+1}}^{2}\geq\frac{1}{64}\sum_{r\geq1}\frac{1}{4^{r}}V_{4^{r}}^{2},
\]
which implies
\[
\sum_{r\geq1}\frac{1}{2^{2r}}V_{2^{2r}}^{2}\leq64\sum_{n\geq1}\frac{V_{n}^{2}%
}{n^{2}}.
\]
Then, by the subadditivity property, we have $V_{2^{2r+1}}\leq2V_{2^{2r}},$ so
that
\[
\sum_{r\geq1}\frac{1}{2^{2r+1}}V_{2^{2r+1}}^{2}\leq\sum_{r\geq1}\frac
{1}{2^{2r}}V_{2^{2r}}^{2}%
\]
and, as a consequence
\[
\sum_{r\geq1}{\frac{1}{2^{2r}}}V_{2^{2r}}^{2}+\sum_{r\geq1}{\frac{1}{2^{2r+1}%
}}V_{2^{2r+1}}^{2}\leq65\sum_{n\geq1}\frac{V_{n}^{2}}{n^{2}}%
\]
and the proof is complete. $\square$

\bigskip

Next lemma contains examples of subadditive sequences which are relevant for
the proofs.

\begin{lemma}
\label{Lsubad} For any stationary Markov chain the sequences $(||E(S_{n}%
|\xi_{0})||)_{n\geq0},$ $(||E(S_{n}|\xi_{n})||)_{n\geq0}$ and $(||E(S_{n}%
|\xi_{0},\xi_{n})||)_{n\geq0}$ are all subadditive.
\end{lemma}

Proof. The proofs are similar. So it is enough to sketch the proof of only one
of them. By the triangle inequality, Markov property and the properties of
conditional expectation, for all positive integers $m$ and $n$,%
\begin{gather*}
||E(S_{n+m}|\xi_{0},\xi_{n+m})||\leq||E(S_{n}|\xi_{0},\xi_{n+m})||+||E(S_{n+m}%
-S_{n}|\xi_{0},\xi_{n+m})||\\
=||E(S_{n}|\mathcal{F}_{0},\mathcal{F}^{n+m})||+||E(S_{n+m}-S_{n}%
|\mathcal{F}_{0},\mathcal{F}^{n+m})||\\
\leq||E(S_{n}|\mathcal{F}_{0},\mathcal{F}^{n})||+||E(S_{n+m}-S_{n}%
|\mathcal{F}_{n},\mathcal{F}^{n+m})||
\end{gather*}
and the result follows. $\square$

\bigskip

Finally we give a lemma for sequences of real numbers.

\begin{lemma}
\label{aux} Let $(a_{i})_{i\geq1}$ be a sequence of real numbers. Then, for
any $m\geq1$%
\[
A_{m}:=\sum\nolimits_{k\geq1}^{m}\frac{1}{k^{2}}(\sum\nolimits_{i=1}^{k}%
a_{i})^{2}\leq4\sum\nolimits_{i=1}^{m}a_{i}^{2}.
\]

\end{lemma}

Proof. Note that%
\begin{align*}
(\sum\nolimits_{i=1}^{k}a_{i})^{2}  &  =\sum\nolimits_{i=1}^{k}a_{i}^{2}%
+2\sum\nolimits_{i=1}^{k}a_{i}\sum\nolimits_{j=1}^{i-1}a_{j}\leq\\
&  2\sum\nolimits_{i=1}^{k}a_{i}\sum\nolimits_{j=1}^{i}a_{j}.
\end{align*}
By changing the order of summation, taking into account that for $i\geq1$ we
have $\sum\nolimits_{k\geq i}k^{-2}\leq2i^{-1}$ and applying the
Cauchy-Schwartz inequality, we obtain
\begin{gather*}
\sum\nolimits_{k=1}^{m}\frac{1}{k^{2}}\sum\nolimits_{i=1}^{k}a_{i}%
\sum\nolimits_{j=1}^{i}a_{j}=\sum\nolimits_{i=1}^{m}a_{i}\sum\nolimits_{j=1}%
^{i}a_{j}\sum\nolimits_{k\geq i}\frac{1}{k^{2}}\\
\leq2\sum\nolimits_{i=1}^{m}|a_{i}|\frac{1}{i}|\sum\nolimits_{j=1}^{i}%
a_{j}|\leq2\left(  \sum\nolimits_{i=1}^{m}a_{i}^{2}\right)  ^{1/2}\left(
\sum\nolimits_{i=1}^{m}\frac{1}{i^{2}}|\sum\nolimits_{j=1}^{i}a_{j}%
|^{2}\right)  ^{1/2}.
\end{gather*}
So%
\[
A_{m}\leq2\left(  \sum\nolimits_{i=1}^{m}a_{i}^{2}\right)  ^{1/2}A_{m}^{1/2},
\]
which completes the proof of the lemma. $\ \square$

\bigskip

\textbf{Acknowledgement.} This paper was partially supported by the NSF grant
DMS-1811373. The author would like to thank Christophe Cuny for discussions
about the square root condition and the variance of partial sums. Many thanks
are going to Michael Lin for pointing out additional references and remarks.
The author would also like to thank the referee for carefully reading the
manuscript and for several observations which improved the presentation of
this paper.

\bigskip$\square$

\end{document}